\documentclass{amsart}

\usepackage[T1]{fontenc}
\usepackage{amsfonts}
\usepackage{amssymb}

\newtheorem{thm}{theorem}[section]
\newtheorem{theorem}[thm]{Theorem}
\newtheorem{proposition}[thm]{Proposition}
\newtheorem{lemma}[thm]{Lemma}
\newtheorem{corollary}[thm]{Corollary}
\newtheorem{remark}[thm]{Remark}

\newtheorem{definition}[thm]{Definition}

\begin{document}

\title[Primeness property for central polynomials of verbally prime algebras]{Primeness property for graded central polynomials of verbally prime algebras}

\author[Silva]{Diogo Diniz}\address{Departamento de Matem\'atica\\ UAME/CCT-UFCG \\ Avenida Apr\'igio Veloso 882\\ 58109-970 Campina Grande-PB, Brazil} \email{diogo@dme.ufcg.edu.br} \thanks{The first author was partially supported by CNPq}

\author[Bezerra]{Claudemir Fidelis Bezerra J\'unior}\address{Departamento de Matem\'atica\\ UAME/CCT-UFCG \\ Avenida Apr\'igio Veloso 882\\ 58109-970 Campina Grande-PB, Brazil} \email{claudemir@mat.ufcg.edu.br}



\begin{abstract}
Let $F$ be an infinite field. The primeness property for central polynomials of $M_n(F)$ was established by A. Regev, i.e., if the product of two polynomials in distinct variables is central then each factor is also central. In this paper we consider the analogous property for $M_n(F)$ and determine, within the elementary gradings with commutative neutral component, the ones that satisfy this property, namely the crossed product gradings. Next we consider $M_n(R)$, where $R$ admits a regular grading, with a grading such that $M_n(F)$ is a homogeneous subalgebra and provide sufficient conditions - satisfied by $M_n(E)$ with the trivial grading - to prove that $M_n(R)$ has the primeness property if $M_n(F)$ does. We also prove that the algebras $M_{a,b}(E)$ satisfy this property for ordinary central polynomials. Hence we conclude that, over a field of characteristic zero, every verbally prime algebra has the primeness property.
\end{abstract}

\maketitle

\section{Introduction}

The study of central polynomials is an important part of the theory of algebras with polynomial identites. Verbally prime algebras were introduced by A. Kemer \cite{Kemer} in his solution of the Specht problem and are of fundamental importance in the theory of p.i. algebras. The existence of central polynomials for verbally prime algebras was proved by Yu. P. Razmyslov \cite{Raz} and earlier for matrix algebras, independently, by Formanek \cite{Formanek} and Razmyslov \cite{Razmyslov}. A p.i. algebra $A$ is verbally prime if whenever $f(x_1,\dots, x_r)$ and $g(x_{r+1},\dots, x_s)$ are two polynomials in distinct variables and $f\cdot g$ is an identity for $A$ then either $f$ or $g$ is an identity for $A$. In his structure theory of $T$-ideals Kemer proves that over a field of characteristic zero every non-trivial verbally prime p.i. algebra is equivalent to one of the algebras $M_n(F)$, $M_n(E)$ (here $E$ denotes the Grassmann algebra of a vector space of countable dimension) and certain subalgebras $M_{a,b}(E)$ of $M_{a+b}(E)$ (see Section \ref{s1}).

Amitsur proved in \cite{Amitsur} that if the field $F$ is infinite the ideal of identities of $M_n(F)$ is a prime ideal. As an analog of Amitsur's result for central polynomials Regev proved in \cite{Regev} the following primeness property on the central polynomials: if $f(x_1,\dots, x_r)$ and $g(x_{r+1},\dots, x_s)$ are two polynomials in distinct variables and $f\cdot g$ is a central polynomial for $M_n(F)$ then both $f$ and $g$ are central.

Graded polynomial identities and graded central polynomials play an important role in the study of p.i. algebras, they were used in the theory developed by Kemer. Such identities are easier to describe in many important cases and they are related to the ordinary ones, for example, coincidence of graded identities implies the coincidence of the ordinary ones. In \cite{D} regular gradings on the algebras $M_n(E)$ and $M_{a,a}(E)$ (see \cite{Regev&Bah}) were used to prove that, if the field $F$ is of characteristic zero, these algebras satisfy the primeness property. In this paper the canonical $\mathbb{Z}_2$-grading on $M_{a,b}(E)$ is used to prove (see Theorem \ref{main}) that these algebras also satisfy the primeness property. This proves that, over a field of characteristc zero, every verbally prime algebra satisfies the primeness property.

We remark that this property arises implicitly in the description of (graded) central polynomials of some verbally prime algebras. For example, over an  infinite field of characteristic $p>2$ the central polynomials of $E$ together with its identities form a so called limit $T$-space and the primeness property is useful to describe these central polynomials (see \cite{BKKS}).

Graded simple algebras, over an algebraically closed field, have been described in \cite{BSZ}, for a group $G$ these are essentially matrix algebras over a twisted group algebra $F^{\alpha}H$, where $\alpha$ is a 2-cocicle with entries in $F^{\times}$. If $H$ is abelian then the canonical $H$-grading on $F^{\alpha}H$ is a regular grading. Classification of abelian gradings on matrix algebras plays an important role in the cassification of gradings on Lie algebras (see \cite{EK}).  A grading on $M_n(F)$ is elementary if every elementary matrix is homogeneous, these gradings play an important role in the description of graded simple algebras. The graded identities and central polynomials of matrix algebras have been described for various gradings. For example, with the canonical $\mathbb{Z}$ and $\mathbb{Z}_n$, central polynomials were described in \cite{B} and more generally in \cite{AK} with crossed product gradings. It is then a natural question to consider the analogous primeness property for the graded polynomials of $M_n(F)$. Graded identities for the algebra $M_n(F)$ with an elementary grading by the group $G$ such that the component associated to the neutral element of $G$ is a commutative subalgebra were described in \cite{BD} and \cite{D}. In Theorem \ref{maing} we prove that with such a grading the algebra $M_n(F)$ satisfies the primeness property if and only if the grading is a crossed product grading (i.e. $G$ has order $n$) and there exists no nontrivial homomorphism $G\rightarrow F^{\times}$.

In Section \ref{ls} we consider matrix algebras $M_n(R)$ with entries in an algebra $R$ that has a regular grading by an abelian group $H$. The gradings considered on $M_n(R)$ are ones for which $R$ and $M_n(F)$ are homogeneous subalgebras and $R$ lies in the neutral component. Under natural hypothesis on the regular grading on $R$ we prove that $M_n(R)$ satisfies the property whenever $M_n(F)$ does. As a corollary we prove that $M_n(E)$ (see Corollary \ref{MnE}) satisfies the primeness property. Finally in Section \ref{smab} the canonical $\mathbb{Z}_2$-grading on $M_{a,b}$ is used to we prove the property for the ordinary central polynomials of these algebras.

\section{Definitions and Preliminary Results}\label{s1}

Let $F$ be an infinite field of characteristic different from $2$, $G$ a group with identity element $\epsilon$ and $A$ an $F$-algebra. A \textit{grading by the group $G$} on $A$ is a vector space decomposition $A=\oplus_{g\in G} A_g$ such that $A_{g}A_{h}\subseteq A_{gh}$ for every $g,h$ in $G$. The subspaces $A_{g}$ are called the \textit{homogeneous components} of $A$ and a non-zero element of $A_{g}$ is said to be homogeneous \textit{of degree} $g$. A subalgebra $B$ is said to be a \textit{homogeneous subalgebra} if $B=\oplus_{g\in G} A_{g}\cap B$. We denote by $F\langle X_G \rangle$ the free $G$-graded algebra, freely generated by the set $X_{G}=\{x_{i,g}|i\in \mathbb{N}, g\in G\}$. This algebra has a natural grading by $G$ where the homogeneous component $(F\langle X_G \rangle)_g$ is the subspace generated by the monomials $x_{i_1,g_1}\cdots x_{i_m,g_m}$ such that $g_{1}\cdots g_m=g$. Let $f(x_{1,g_1},\dots, x_{m,g_m})$ be a polynomial in $F\langle X_G \rangle$, an $m$-tuple $(a_1,\dots, a_m)$ such that $a_i\in A_{g_i}$ for $i=1,\dots,m$ is called an \textit{$f$-admissible substitution} (or simply an admissible substitution). If $f(a_1,\dots,a_m)=0$ for every admissible substitution $(a_1,\dots, a_m)$ we say that the polynomial $f$ is a \textit{graded polynomial identity} for $A$.  We denote by $Z(A)$ the centre of $A$ and we say that $f$ is a \textit{graded central polynomial} for $A$ if $f$ is not a graded polynomial identity and $f(a_1,\dots, a_m)\in Z(A)$ for every admissible substitution $(a_1,\dots, a_m)$.

\begin{definition}
Let $A=\oplus_{g\in G}A_g$ be an algebra with a grading by the group $G$. We say that $A$ satisfies the \textbf{primeness property} for (graded) central polynomials if, whenever the product $f(x_{1,g_1},\cdots, x_{r,g_{r}})\cdot g(x_{r+1,g_{r+1}},\cdots, x_{s,g_{s}})$ of two polynomials in distinct variables is a graded central polynomial, then both $f(x_{1,g_1},\cdots, x_{r,g_{r}})$ and $ g(x_{r+1,g_{r+1}},\cdots, x_{s,g_{s}})$ are graded central polynomials for $A$.
\end{definition}

The set of all graded central polynomials and all graded identities for $A$ is a subspace of $F\langle X_G \rangle$ invariant under all endomorphisms of $F\langle X_G \rangle$, such subspaces are called \textit{$T_G$-spaces}. The intersection of $T_G$-spaces is also a $T_G$-space. Given a subset
$S\subset F\langle X_G \rangle$ we denote by $\langle S \rangle_{T_G}$ the intersection of all $T_G$-spaces that contain $S$, this is the \textit{$T_G$-space generated by $S$}. If $S=\{f\}$ we use the notation $\langle f \rangle_{T_G}$ for the $T_G$-space generated by $\{f\}$. We remark that  $\langle f \rangle_{T_G}$ is generated as a vector space by the set $\{f(p_1,\dots,p_n)|p_i \in (F\langle X_G \rangle)_{g_i}\}$. If $f$ is a graded central polynomial for $A$ then every element of $\langle f \rangle_{T_G}$ either a graded identity or a graded central polynomial for $A$. If $G=\{\epsilon\}$ we recover the definition of ordinary polynomial identities and central polynomials, in this case we use the notation $F\langle X \rangle$ for the free associative algebra and $x_{i}$ for the variables.

The Grassmann algebra $E$ of a countable dimensional vector space with basis $\{e_1, e_2,\dots\}$ has the presentation:
\[\langle 1, e_1, e_2, \dots| e_ie_j=-e_je_i,\mbox{ for all i,j} \geq 1 \rangle. \] It is well known that $\mathcal{B}_E=\{1,e_{i_1}e_{i_2}\dots e_{i_k}|i_1<i_2<\dots <i_k\}$ is a basis for $E$. Hence we have $E=E_{0}\oplus E_{1}$ where $E_{0}$ is the subspace generated by $1$ and all monomials with even $k$ and $E_1$ is the subspace generated by the monomials with odd $k$. This decomposition is a $\mathbb{Z}_2$-grading on $E$. The algebra $M_{a+b}(E)$ has a subalgebra denoted $M_{a,b}(E)$ which consists of matrices of the form
\[
\left(\begin{array}{cccc}
	A & B\\
	C & D
\end{array}\right),
\] where $A\in M_{a}(E_0)$, $D \in M_b(E_0)$ and $B$, $C$ are blocks with entries in $E_1$. This algebra has a canonical $\mathbb{Z}_2$-grading with homogeneous components $\left(M_{a,b}(E)\right)_0$ and $\left(M_{a,b}(E)\right)_1$ that consists of matrices of the form \[
\left(\begin{array}{cccc}
	A & 0\\
	0 & D
\end{array}\right) \mbox{\hspace{0.6cm} and \hspace{0.6cm}} \left(\begin{array}{cccc}
	0 & B\\
	C & 0
\end{array}\right),
\] respectively.

\begin{proposition}\label{peq}
Let $f(x_1,\dots, x_r)$ be a multihomogeneous polynomial of degree $d_i$ in $x_i$ and $h_i(x_1,\dots, x_r, z_i)$ the multihomogeneous component of $f(x_1+z_1,\dots, x_r+z_r)$ of degree one in $z_i$. Then
\begin{equation}
[f(x_1,\dots, x_r),x_{r+1}]=\sum_i h_i(x_1,\dots, x_r,[x_i,x_{r+1}]).
\end{equation}
\end{proposition}

\textit{Proof.}
The equality follows from the fact that the map $a\mapsto [a,x_{r+1}]$ is a derivation of $F\langle X \rangle$ (see[pg. 9]\cite{GZ}).
\hfill $\Box$

\begin{proposition}\label{C}
Let $A$ be a $G$-graded algebra with basis $\mathcal{B}$ that consists of homogeneous elements, $V$ a subspace of $A$ and $f(x_{1,g_1},\dots, x_{n,g_n})$ a polynomial. The following assertions are equivalent:
\begin{enumerate}
\item[(i)] For every admissible substitution $a_1,\dots, a_n$ by elements of $A$, the element $f(a_1,\dots, a_n)$ lies in $V$;
\item[(ii)] For every polynomial $f^{\prime}(x_{1,g_1},\dots, x_{m, g_m})$ in $\langle f \rangle_{T_G}$ and every admissible substitution $a_1,\dots, a_m$ by elements of $A$, the element $f^{\prime}(a_1, \dots, a_m)$ lies in $V$.
\item[(iii)] For every multihomogeneous polynomial $f^{\prime}(x_{1,g_1},\dots, x_{m, g_m})$ in $\langle f \rangle_{T_G}$ and every admissible substitution $b_1,\dots, b_m$ by elements of the basis $\mathcal{B}$, the element $f^{\prime}(b_1, \dots, b_m)$ lies in $V$.
\end{enumerate}
\end{proposition}
\textit{Proof}
The set $S=\{f(p_1,\dots, p_n)|p_i \in (F\langle X_G \rangle)_{g_i}\}$ generates $\langle f \rangle_{T_G}$ as a vector space. If the assertion $(i)$ holds for $f$ then it holds for every element of $S$ and therefore for every $f^{\prime}$ in $\langle f \rangle_{T_G}$. Hence we conclude $(ii)$. It is clear that $(ii)$ implies $(iii)$. Now we assume that $(iii)$ holds. Let $a_1,\dots, a_n$ an admissible substitution by elements of $A$. We write $a_i=\sum_{j=1}^{n_i} \lambda_{ij}b_j$, $i=1,\dots, n$, where $b_j\in \mathcal{B}$. Let $s_1=0$ and $s_{i+1}=s_{i}+n_i$ for $i=2,\dots, m-1$. Denote $q_i(x_{s_i+1,g_i},\dots, x_{s_i+n_i,g_i})=\sum_{j=1}^{n_i}\lambda_{ij}x_{s_i+j, g_i}$ and write \[f(q_1,\dots q_m)=f_1+\dots+f_l,\] as a sum of multihomogeneous components. The field $F$ is infinite and therefore the polynomials $f_1,\dots, f_l$ lie in $\langle f \rangle_{T_G}$. Let $S$ be the substitution such that $x_{s_i+j, g_i}$ is replaced by $b_j$. The result $q_j|_S$ of this substitution is $q_j(b_{1},\dots, b_{n_i})=a_i$ and therefore $f(a_1,\dots, a_n)=f_1|_S+\dots+ f_l|_S$. From $(iii)$ we conclude that the elements $f_1|_S,\dots, f_l|_S$ lie in $V$, hence $f(a_1,\dots, a_n)$ also lies in $V$.
\hfill $\Box$

\begin{corollary}\label{CV}
Let $A$ be an algebra graded by the group $G$, $V$ a subspace of $A$ and let $f(x_{1,g_1},\dots, x_{m,g_m})$ be a polynomial such that for every admissible substitution $a_1,\dots, a_m$ by elements of $A$, the element $f(a_1,\dots, a_m)$ lies in $V$. Then for every $b$ in $A_{\epsilon}$ and every admissible substitution $a_1,\dots, a_m$ the commutator $[f(a_1,\dots, a_m),b]$ also lies in $V$.
\end{corollary}

\textit{Proof.}
Proposition \ref{peq} implies that
\begin{equation}\label{e22g}
[f(x_{1,g_1},\dots, x_{m,g_m}),x_{m+1,\epsilon}]=\sum_i h_i(x_{1,g_1},\dots, x_{m,g_m},[x_{i,g_i},x_{m+1,\epsilon}]).
\end{equation}
The variable $x_{m+1,\epsilon}$ is homogeneous of degree $\epsilon$, therefore $[x_{i,g_i},x_{m+1,\epsilon}]$ has degree $g_i$ and the polynomial $h_i(x_{1,g_1},\dots, x_{m,g_m},[x_{i,g_i},x_{m+1,\epsilon}])$ lies in the $T_G$-space generated by $f$. The substitution $(a_1,\dots, a_m,b)$ is $h_i$-admissible, hence the previous proposition implies that $h_i(a_1,\dots, a_m,[a_i,b])$ lies in $V$. The result now follows directly from (\ref{e22g}).
\hfill $\Box$

\begin{proposition}\label{ll}
Let $f(x_{1, g_1},\dots, x_{r,g_r})$, $g(x_{r+1, g_{r+1}},\dots, x_{s,g_s})$ be polynomials in distinct variables. If $f^{\prime}\in \langle f \rangle_{T_G}$ and $g^{\prime}\in \langle g \rangle_{T_G}$ then $f^{\prime}g^{\prime}$ lies in $\langle f\cdot g \rangle_{T_G}$. In particular if $f\cdot g$ is a central polynomial for $A$ then $f^{\prime}g^{\prime}$ is also a central polynomial for $A$.
\end{proposition}

\textit{Proof.}
	Note that $f^{\prime}$ and $g^{\prime}$ are linear combinations of polynomials of the form $f(p_{1},\dots, p_{r})$ and $g(q_{1},\dots, q_{s})$ respectively, where $p_1,\dots, p_r$ and $q_1,\dots,q_s$ admissible substitutions by polynomials in $F\langle X_G \rangle$. Hence the product $f^{\prime}\cdot g^{\prime}$ is a linear combination polynomials in the set $$S=\{f(p_{1},\dots, p_{r})\cdot g(q_{1},\dots, q_{r})|p_1,\dots, p_r, q_1,\dots, q_s \in F\langle X \rangle\}.$$ Since $f$ and $g$ are polynomials in distinct variables each element in $S$ lies in $\langle f\cdot g \rangle_T$. Therefore $f^{\prime}\cdot g^{\prime}$ lies in $\langle f\cdot g \rangle_T$. The second assertion is an obvious consequence of the first one.
\hfill $\Box$

\section{Elementary gradings on $M_n(F)$ and the primeness property}

In this section $A$ denotes the algebra $M_n(F)$. We denote $E_{ij}$ the elementary matrix with $1$ in the $(i,j)$-th entry and $0$ in the other entries.

\begin{definition}
Let $\textbf{g}=(g_1,\dots, g_n)$ be an $n$-tuple of elements in the group $G$. The decomposition $A=\oplus_{g\in G} A_g$, where $A_g$ is the subspace generated by the elementary matrices $E_{ij}$ such that $g_{i}^{-1}g_j=g$, is a grading by $G$ on $A$ called the \textbf{elementary grading }induced by $\textbf{g}$.
\end{definition}

We recall that an elementary grading by an $n$-tuple of parwise distinct elements by a group $G$ of order $n$ is called a \textit{crossed product grading} (see \cite{AK}). Our main result in this section characterizes the elementary gradings by an $n$-tuple of pairwise distinct elements on $A$ for which this algebra satisfies the primeness property for graded central polynomials. We prove in Theorem \ref{maing} that these gradings are crossed product gradings by a group $G$ that has no non-trivial homomorphism $G\rightarrow F^{\times}$.

\begin{proposition}\label{pp}
Let $A$ be the matrix algebra $M_n(F)$ with an elementary grading and $f(x_{1, g_1},\dots, x_{r, g_r})$, $g(x_{r+1, g_{r+1}},\dots, x_{s, g_s})$ polynomials in distinct variables such that the product $f\cdot g$ is a graded central polynomial for $A$. Then there exists an invertible diagonal matrix $P$ such that the result of every $f$-admissible (resp. $g$-admissible) substitution is a scalar multiple of $P$ (resp. $P^{-1}$).
\end{proposition}
\textit{Proof.}
Let $f(x_{1,g_1},\dots, x_{r,g_r})$ and $g(x_{r+1,g_{r+1}},\dots, x_{s,g_s})$ be polynomials in distinct variables such that $f\cdot g$ is a graded central polynomial for $A$. There exists an admissible substitution $A_1,\dots, A_s$ such that $f(A_1,\dots, A_r)\cdot g(A_{r+1},\dots, A_s)$ is a non-zero scalar matrix. Hence $g(A_{r+1},\dots, A_s)$ is an invertible matrix, let $P$ denote its inverse. If $B_1,\dots, B_s$ is a $f$-admissible substitution the product $f(B_1,\dots, B_r)\cdot g(A_{r+1},\dots, A_s)$ is a scalar matrix. Therefore $f(B_1,\dots, B_r)$ is a scalar multiple of $P$ for every $f$-admissible substitution. Analogously, the result of every $g$-admissible substitution is a scalar multiple of $P^{-1}$. We now prove that $P$ is a diagonal matrix.

It follows from Corollary \ref{CV} that if $M$ is any homogeneous matrix of degree $\epsilon$ and $B_1,\dots, B_r$ is an $f$-admissible substitution  the commutator $[f(B_1,\dots, B_r), M]$ is a scalar multiple of $P$. Therefore $[P, M]=\lambda_M P$ for some scalar $\lambda_M$. The grading is elementary, hence the matrices $E_{ii}$, $i=1,\dots, n$ have degree $\epsilon$ and for each $i$ there exists a scalar $\lambda_i$ such that
\begin{equation}\label{e0}
[P,E_{ii}]=\lambda_{i}P.
\end{equation}
We write $P=\sum_{i,j}p_{ij}E_{ij}$. If $\lambda_i=0$ for every $i$ then an easy computation yields that $P$ is a diagonal matrix. We prove that this must be the case. Assume that $\lambda_i\neq 0$ for some $i$. The left side of (\ref{e0}) is $\sum_{k\neq i}p_{ki}E_{ki}-\sum_{l\neq i}p_{il}E_{il}$. Since $\lambda_i\neq 0$ we conclude that $p_{ii}=0$ and $p_{kl}=0$ if $k\neq i$ and $l\neq i$. Moreover $p_{ki}=\lambda_ip_{ki}$ for every $k\neq i$ and $-p_{il}=\lambda_i p_{il}$ for every $l\neq i$. If $p_{ki}\neq 0$ for some $k\neq i$ then $\lambda_i=1$ and $p_{il}=0$ for every $l\neq i$. In this case $P=\sum_{k\neq i}p_{ki}E_{ki}$ which is a contradition since this matrix is not invertible. If $p_{il}\neq 0$ for some $l\neq i$ then $\lambda_i=-1$ and $P=\sum_{l\neq i}p_{il}E_{il}$, this last matrix is not invertible and again we have a contradiction. The remaining possibility is that $p_{ki}=0$ for every $k\neq i$ and $p_{li}=0$ for every $l\neq i$ which implies $P=0$. This is also a contradiction.
\hfill $\Box$

Our main objective now is to prove that for elementary gradings determined by tuples of pairwise distinct elements the primeness property holds only for crossed product gradings by suitable groups. The precise statement is in Theorem \ref{maing}.

\begin{definition}
Let $A$ be the matrix algebra $M_n(F)$ with an elementary grading. We denote $\emph{H}$ the set of permutations $\sigma$ in $S_n$ for which the corresponding automorphism $\Lambda_{\sigma}$, given by $\Lambda_{\sigma}(E_{ij})=E_{\sigma(i)\sigma(j)}$, is a graded automorphism of $A$.
\end{definition}

The map $\sigma \mapsto \Lambda_{\sigma}$ is a homomorphism and the set $\emph{H}$ in the definition above is a subgroup of $S_n$. Next we prove some results that involve this subgroup.

\begin{proposition}\label{r}
If $f$ is a non-identity such that the result of every $f$-admissible substitution is a scalar multiple of $P$, it follows that for every $\sigma \in \emph{H}$ there exists a non-zero scalar $\lambda(\sigma)$ such that \[\Lambda_\sigma (P)= \lambda(\sigma)^{-1} P. \] The map $\lambda:\emph{H}\rightarrow F^{\times}$ is a homomorphism.
\end{proposition}
\textit{Proof.}
There exists an admissible substitution $f(A_1,\dots, A_m)=\lambda P$ for some scalar $\lambda \neq 0$. The automorphisms $\Lambda_{\sigma}$ are graded and therefore \[ \Lambda_{\sigma}(\lambda P)=\Lambda_{\sigma}(f(A_1,\dots, A_m))=f(\Lambda_{\sigma}(A_1),\dots, \Lambda_{\sigma}(A_m))= \lambda^{\prime}P,\] for some non-zero $\lambda^{\prime}$ depending on $\sigma$. Hence there exists a non-zero scalar $\lambda_{\sigma}$ such that $\Lambda_{\sigma}(P)=\lambda_{\sigma}^{-1}P$. It is clear that $\lambda:\emph{H}\rightarrow F^{\times}$ given by $\sigma\mapsto \lambda_{\sigma}$ is a homomorphism.
\hfill $\Box$

In the next lemma we describe the matrices $P$ that satisfy the equality $\Lambda_\sigma (P)= \lambda(\sigma)^{-1} P$ in terms of the homomorphism $\lambda$ and of the canonical action of $\emph{H}$ on $\{1,2,\dots,n\}$.

\begin{lemma}\label{l2}
Let $\lambda:\emph{H}\rightarrow F^{\times}$ be a homomorphism and for each $i$ denote by $P(i)=\sum_{\alpha\in \emph{H}}\lambda(\alpha)E_{\alpha(i), \alpha(i)}$. The following assertions hold:
\begin{enumerate}
\item[(i)] For every $\sigma \in \emph{H}$ we have $\Lambda_{\sigma}(P(i))=\lambda(\sigma)^{-1}P(i)$;
\item[(ii)] The matrices $P(i)$ and $P(j)$ have a non-zero entry in the same position if and only if $i$ and $j$ lie in the same $\emph{H}$-orbit;
\item[(iii)] Let $i_1,\dots, i_d$ be a system of representatives of the orbits of the canonical action of $H$ on $\{1,2,\dots, n\}$. The diagonal matrix $P$ satisfies
\begin{equation}\label{sca}
\Lambda_\sigma (P)= \lambda(\sigma)^{-1} P
\end{equation}
for every $\sigma \in \emph{H}$ if and only if $P$ is a linear combination of the matrices $P(i_1),\dots, P(i_d)$.
\end{enumerate}
\end{lemma}
\textit{Proof.}
We have \[\Lambda_{\sigma}(P(i))=\sum_{\alpha\in \emph{H}}\lambda(\alpha)E_{\sigma\alpha(i), \sigma\alpha(i)}=\lambda(\sigma)^{-1}\sum_{\alpha\in \emph{H}}\lambda(\sigma\alpha)E_{\sigma\alpha(i), \sigma\alpha(i)}=\lambda(\sigma)^{-1}P(i),\] hence assertion $(i)$ holds. Assertion $(ii)$ is clear since the non-zero entries of $P(i)$ are in the positions corresponding to the orbit of $i$. Next we prove $(iii)$. Let $P=\sum_i p_i E_{ii}$ be a matrix that satisfies (\ref{sca}) and $\lambda_1,\dots, \lambda_d$ be scalars such that the $i_s$-th entry of the matrix $\lambda_s P(i_s)$ is $p_{i_s}$. Let $P^{\prime}=\sum_s \lambda_s P(i_s)$ and write $P^{\prime}=\sum_i p_i^{\prime}E_{ii}$. It follows from $(ii)$ that $p_{i_s}=p_{i_s}^{\prime}$ for $s=1,\dots, d$. Assertion $(i)$ implies that $P^{\prime}$ also satisfies (\ref{sca}), hence  $p_{\sigma(i)}=p^{\prime}_{\sigma(i)}$ whenever $p_i=p_i^{\prime}$ and $\sigma \in \emph{H}$. Since $i_1,\dots, i_d$ is a system of representatives this implies that $P=P^{\prime}$.
\hfill $\Box$

For elementary gradings by tuples of pairwise distinct elements it is possible to provide for the matrices in the above lemma a non-identity $f$ such that the result of every admissible substitution is a scalar multiple of $P$. This is done in the next proposition.

\begin{lemma}\label{subs}
Let $x_{1,g_1}\cdots x_{n,g_n}$ be a multilinear monomial of degree $\epsilon$ such that the $n$-tuple $(g_1, g_1g_2,\dots, g_1\cdots g_n)$ consists of pairwise distinct elements.  If the $n$-tuples $(E_{i_1,j_1},\dots, E_{i_n,j_n})$ and $(E_{k_1,l_1},\dots, E_{k_n,l_n})$ are admissible substitutions and the result of each substitution is not zero then there exists $\sigma \in \emph{H}$ such that \[(E_{k_1,l_1},\dots, E_{k_n,l_n})=(E_{\sigma(i_1),\sigma(j_1)},\dots, E_{\sigma(i_n),\sigma(j_n)}).\]
\end{lemma}
\textit{Proof.}
Since  $E_{i_1,j_1}\cdots E_{i_n,j_n}\neq 0$ and the total degree of the polynomial is $\epsilon$ and we conclude that $j_1=i_2, \dots, j_{n-1}=i_n, j_n=i_1$. Moreover the fact that $(g_1, g_1g_2,\dots, g_1\cdots g_n)$ consists of pairwise distinct elements implies that $(i_1,\dots, i_n)$ is a permutation of the elements $1,\dots, n$. Analogously $(k_1,\dots, k_n)$ is a permutation. We denote by $\sigma$ the permutation in $S_n$ such that $\sigma(i_s)=k_s$, for $s=1,\dots, n$. The substitutions considered are admissible, therefore the automorphism of $A$ such that $E_{ij}\mapsto E_{\sigma(i),\sigma(j)}$ maps the matrices $E_{i_1,i_2},\dots, E_{i_{n-1}, i_{n}}, E_{i_n,i_1}$ onto elementary matrices of the same degree. This implies that this automorphism is a graded automorphism.
\hfill $\Box$

\begin{proposition}\label{p1}
Let $A$ be algebra $M_n(F)$ with an elementary grading induced by the $n$-tuple $(h_1,\dots, h_n)$ of pairwise distinct elements. If there exists a non-zero diagonal matrix $P$ and a non-identity $f(x_{1,g_1},\dots, x_{m,g_m})$ such that the result of every admissible substitution is a scalar multiple of $P$ then there exists a homomorphism $\lambda:\emph{H}\rightarrow F^{\times}$ such that $P$ equals a linear combination of the matrices $P(i_1),\dots, P(i_d)$ in the previous lemma. Conversely let $\lambda:\emph{H}\rightarrow F^{\times}$ be a homomorphism and $P=diag (p_1,\dots, p_n)$ a non-zero linear combination of the matrices $P(i_1),\dots, P(i_d)$, then the polynomial
\begin{equation}\label{polynomial}
f(x_{1,g_1},\dots, x_{n,g_n})=\sum_{i} p_{i} x_{\sigma^{i}(1),g_{\sigma ^{i}(1)}},\dots, x_{\sigma ^{i}(n),g_{\sigma ^{i}(n)}},
\end{equation}
where $\sigma$ is the $n$-cicle $(1 2 \dots  n)$, $g_i=h_i^{-1}h_{i+1}$ for $i=1,\dots, n-1$ and $g_n=h_n^{-1}h_1$; is not an identity and has the property that every the result of every admissible substitution is a scalar multiple of $P$.
\end{proposition}

\textit{Proof.}
The first claim about $P$ follows form Proposition \ref{r} and $(iii)$ of the Lemma \ref{l2}.
Now we prove the second claim. The polynomial $f$ is multilinear, hence it suffices to prove that the result of every admissible substitution by elementary matrices is a scalar multiple of $P$. It is clear that $f(E_{12},\dots,E_{n1})=P$, in particular $f$ is not a graded identity. Lemma \ref{subs} implies that every admissible substitution such that the result is not zero is of the form $(\Lambda_{\sigma}(E_{12}),\dots, \Lambda_{\sigma}(E_{n1}))$. Now $(iii)$ of the Lemma \ref{l2} implies \[f(\Lambda_{\sigma}(E_{12}),\dots, \Lambda_{\sigma}(E_{n1}))=\Lambda_{\sigma}(P)=\lambda(\sigma)^{-1}P,\] and the second claim follows.
\hfill $\Box$

\begin{lemma}\label{l1}
If $A$ has an elementary grading induced by an $n$-tuple of pairwise distinct elements of $g$ then the orbit of every element in the set $\{1,2,\dots,n\}$ under the canonical action has $|H|$ elements. In particular $|H|$ divides $n$.
\end{lemma}
\textit{Proof.}
If $\sigma(i)=i$ for some $\sigma\in \emph{H}$ and some $i$ then for every $j$ the matrices $E_{i,j}$ and $E_{i,\sigma(j)}$ have the same degree. Since the $n$-tuple that determines the grading consists of pairwise distinct elements this implies that $\sigma(j)=j$. Hence the stabilizer of every element is trivial and the result follows.
\hfill $\Box$

\begin{proposition}\label{p2}
Let $A$ be algebra $M_n(F)$ with an elementary grading induced by the $n$-tuple $(g_1,\dots, g_n)$ of pairwise distinct elements of the group $G$. Then $A$ has a crossed product grading by the support if and only if $|\emph{H}|=n$, moreover in this case the groups $\emph{H}$ and $G$ are isomorphic.
\end{proposition}
\textit{Proof.}
We assume first that $A$ has a crossed product grading. For each $g\in G$ let $\sigma_g\in S_n$ be the permutation such that $$(gg_1,\dots, gg_n)=(g_{\sigma_g(1)}, \dots, g_{\sigma_g(n)}).$$ The map $g\mapsto \sigma_g$ is a group homomorphism. To verify this let $g,h$ be elements in $G$, we have $g_{\sigma_g(j)}=gg_j$ for every $j$ and therefore \[g_{\sigma_g(\sigma_h (i))}=g g_{\sigma_h (i)}=ghg_i=g_{\sigma_{gh}(i)}.\] Hence $\sigma_g \sigma_h =\sigma_{gh}$. If $\sigma_g(j)= j$ for some $j$ we have $g_j=gg_j$, hence $g=\epsilon$. This proves that the homomorphism is injective. Now we prove that its image lies in $\emph{H}$. The degree of $\Lambda_{\sigma_g}(E_{ij})=E_{\sigma_g(i)\sigma_g(j)}$ is \[g_{\sigma_g(i)}^{-1}g_{\sigma_g(j)}=(gg_i)^{-1}(gg_j)=g_i^{-1}g_j,\] hence $\Lambda_{\sigma_g}$ is a graded automorphism. Lemma \ref{l1} implies that $\emph{H}$ has order at most $n$, therefore $g\mapsto \sigma_g$ is also surjective. Hence $\emph{H}$ and $G$ are isomorphic and $|\emph{H}|=n$.

Now we assume that $|\emph{H}|=n$. Lemma \ref{l1} implies that there is only one orbit, i.e., the action of $\emph{H}$ on $\{1,\dots, n\}$ is transitive. Let $g$ be in the support of the grading and $E_{kl}$ be of degree $g$. For every $i$ there exists $\sigma \in \emph{H}$ such that $\sigma(k)=i$, hence $E_{ij}$ has degree $g$ for $j=\sigma(l)$. Now if $h$ also lies in the support of the grading then there exists a matrix $E_{jt}$ of degree $h$ and $E_{it}=E_{ij}E_{jt}$ has degree $gh$. Therefore $gh$ lies in the support, it is clear that $h^{-1}$ lies in the support if $h$ does, hence the support is a subgroup of $G$ of order $n$. Since the elements $(\epsilon, g_1^{-1}g_2,\dots, g_1^{-1}g_n)$ are pairwise distinct elements in the support we conclude that the grading is a crossed product grading by the support.
\hfill $\Box$

\begin{theorem}\label{maing}
Let $A$ be the matrix algebra $M_n(F)$ with an elementary grading by an $n$-tuple of pairwise distinct elements of $G$. The algebra $A$ has the primeness property for central polynomials if and only if the grading is a crossed product grading and the group $G$ has no non-trivial homomorphism $G\rightarrow F^{\times}$.
\end{theorem}

\textit{Proof.}
If $A$ does not have a crossed product grading then it follows from Proposition \ref{p2} that $|\emph{H}|<n$. In this case Lemma \ref{l1} implies that the number $d$ of orbits in the $\emph{H}$ action on $\{1,2,\dots,n\}$ is greater than $1$. Let $i_1,\dots, i_d$ be a system of representatives for the orbits, we apply Proposition \ref{p1} with the trivial homomorphism to obtain a polynomial $f$ such that the result of every admissible substitution is a scalar multiple of the matrix $P=P(i_1)+\dots+P(i_{d-1})-P(i_d)$. Since $d>1$, the polynomial $f$ is not graded central. Clearly $P^2=I$, therefore the product of two copies of $f$ in distinct variables is central. Hence $A$ does not have the primeness property for central polynomials. Now assume that $A$ has a crossed product grading and that $G$ has a non-trivial homomorphism $G\rightarrow F^{\times}$. Proposition \ref{p2} implies that $\emph{H}$ is isomorphic to $G$. In this case there exists a non-trivial homomorphism $\lambda:\emph{H}\rightarrow F^{\times}$. The action of $\emph{H}$ determines only one orbit, we apply Lemma \ref{l2} with this homomorphism and Proposition \ref{p1} to obtain a matrix $P=P(i_1)$ and a polynomial $f$ such that the image of every admissible substitution is a scalar multiple of $P$. The entries in the diagonal of $P$ are in the image of $\lambda$ and therefore are $n$-th roots of the unit. This implies that $P^{n}=I$, hence the product of $n$ copies of $f$ in distinct variables is central. The homomorphism $\lambda$ is not trivial, hence $P$ is not a scalar matrix and $f$ is not a graded central polynomial. Therefore $A$ does not have the primeness property for central polynomials.

Finally assume that $A$ has a crossed product grading and that the only homomorphism $G\rightarrow F^{\times}$ is the trivial one. In this case we fix two polynomials $f(x_{1, g_1},\dots, x_{r, g_r})$ and $g(x_{r+1, g_{r+1}},\dots, x_{s, g_s})$ in distinct variables such that the product $f\cdot g$ is a graded central polynomial for $A$. If we prove that the matrix $P$ in Proposition \ref{pp} is a scalar matrix it follows that $f$ and $g$ are central polynomials. Proposition \ref{r} implies that there exits a homomorphism $\lambda:\emph{H}\rightarrow F^{\times}$ such that \[\Lambda_\sigma (P)= \lambda(\sigma)^{-1} P. \] Since the grading is a crossed product grading $\emph{H}$ is isomorphic to $G$. The only one dimensional homomorphism of $G$ onto $F^{\times}$ is the trivial one, hence $\lambda$ is the trivial homomorphism. The order of $\emph{H}$ is $n$, hence Lemma \ref{l1} implies that the action of $\emph{H}$ on $\{1,\dots, n\}$ determines only one orbit. In this case Lemma \ref{l2} implies that $P$ is a scalar matrix.
\hfill $\Box$

Let $A=\oplus_{g\in G}A_g$ and $A=\oplus_{h\in H}A_h$ be gradings on $A$ by the groups $G$ and $H$ respectively. We say that the $H$-grading on $A$ is a \textit{coarsening} of the $G$-grading on $A$ if for every $g\in G$ there exists $h\in H$ such that $A_g\subseteq A_h$. The next result proves that the primeness property is inherited by coarsenings of a grading with this property.

\begin{theorem}\label{t2}
Let $A=\oplus_{g\in G}A_g$ be a grading on $A=M_n(F)$ by the group $G$ and $A=\oplus_{h\in H}A_h$ an $H$-grading that is a coarsening of the  $G$-grading. If $A$ has the primeness property for $G$-graded central polynomials then it has the primeness property for $H$-graded central polynomials.
\end{theorem}
\textit{Proof.}
Let $f(x_{1, h_1},\dots, x_{r, h_r})$ and $g(x_{r+1, h_{r+1}},\dots, x_{s, h_s})$ be polynomials in distinct variables such that the product $f\cdot g$ is an $H$-graded central polynomial for $A$. As in the proof of Proposition \ref{pp} we obtain a $g$-admissible substitution such that $g(A_{r+1},\dots, A_s)$ is an invertible matrix. Let $(A_1,\dots, A_r)$ be an $f$-admissible substitution. Since the grading by $H$ is a coarsening the grading by $G$ each matrix $A_i$ is a sum of homogeneous matrices in the $G$-grading. Therefore we may write $A_i=B_{1,i}+\dots+B_{k_i,i}$ as a sum of elements of degree $g_{1,i},\dots, g_{k_i,i}$, here $g_{k,i}\neq g_{l,i}$ if $k\neq l$. We write $p_{i}=\sum_{l=1}^{k_i}{x_{i,g_{l,i}}}$ and let $\overline{f}=f(p_1,\dots, p_r)$, $\overline{g}=g(p_{r+1},\dots, p_{s})$. Notice that $p_i$, $p_j$ are in distinct variables if $i\neq j$. Hence $\overline{f}$ and $\overline{g}$ are in distinct variables. It is clear that if we replace $x_{i,g_{l,i}}$ by $B_{l,i}$ in $p_i$ the result is $A_i$ and therefore the results of this substitution in $\overline{f}$ and $\overline{g}$ are $f(A_1,\dots, A_r)$ and $g(A_{r+1},\dots, A_s)$ respectively. If the product $\overline{f}\cdot \overline{g}$ is a graded identity then $f(A_1,\dots, A_r)\cdot g(A_{r+1},\dots, A_s)=0$ and the hypothesis that $g(A_{r+1},\dots, A_s)$ is invertible implies that $f(A_1,\dots, A_r)=0$. Otherwise $\overline{f}\cdot \overline{g}$ is a central polynomial because $f\cdot g$ is central. The primeness property for $G$-graded central polynomials implies that $\overline{f}$ is central and therefore $f(A_1,\dots, A_r)$ is a central element. This proves that $f$ is an $H$-graded central polynomial. Analogously one proves that $g$ is central.
\hfill $\Box$

\section{Matrix algebras over algebras with a regular grading}\label{ls}

Let $R$ be an algebra (with unit) and $A$ the matrix algebra $M_n(R)$. We identify the scalar matrix $\mathrm{diag}\ (r,\dots, r)$ with $r$ and consider $R$ as a subalgebra of $A$. Notice that the matrix algebra $B=M_n(F)$ is also a subalgebra of $A$.  In this section we consider gradings on $A$ such that $R$ and $B$ are homogeneous subalgebras and $R$ has the trivial grading. In this case it is clear that the grading $A=\oplus_{g\in G}A_g$ on $A$ is such that $A_{g}=\{rb|r\in R, b\in B_g\}$. Our main result of this section (Theorem \ref{MnA}) is that if $R$ has a suitable regular grading then $A$ has the primeness property whenever $B$ does. As a corollary we obtain the primeness property for $M_n(E)$ over infinite fields. Next we recall the definition of a regular grading (see \cite{AO}, \cite{Regev&Seman}).

\begin{definition}\label{def.rel}
Let $R$ be an algebra and \[
R=\oplus_{h\in H} R_h,
\]
a grading by the abelian group $H$. The $H$-grading above on $R$ is said to be regular if there exists a commutation function $\beta:H\times H \rightarrow F^{\times}$ such that
\begin{enumerate}
  \item [(P1)] For every $n$ and every $n$-tuple $(h_1,\dots, h_n)$ of elements of $H$, there exist $r_1,\dots,r_n$ such that $r_j\in R_{h_j}$, and $r_1\cdots r_n\neq0$.
  \item [(P2)] For every $h_1,h_2 \in H$ and for every $a_{h_1}\in A_{h_1}$, $a_{h_2}\in A_{h_2}$, we have
      \[
      a_{h_1}a_{h_2}=\beta(h_1,h_2)a_{h_2}a_{h_1}.
      \]
\end{enumerate}
The regular $H$-grading on $R$ is said to be minimal if for any $h\neq \epsilon$ there exists $h^{\prime}$ such that $\beta(h,h^{\prime})\neq 1$.
\end{definition}

\begin{remark}
In the definition above the commutation function $\beta:H\times H \rightarrow F^{\times}$ is a skew-symmetric bicharacter (see \cite[Remark 13]{AO}), i.e., for every $h_0 \in H$ the maps $h\mapsto \beta(h_0,h)$ and $h\mapsto \beta(h,h_0)$ are characters of the group $H$ and $\beta(h_2,h_1)=\beta(h_1,h_2)^{-1}$ for every $h_1,h_2 \in H$.
\end{remark}

In most of the known examples of a minimal regular grading by an abelian group (see \cite{Regev&Bah}) it turns out that the centre of the algebra coincides with the neutral component. In the next proposition we prove this coincidence under a minor condition on the regular grading.

\begin{proposition}\label{center}
Let $R$ be an algebra with a regular grading by the group $H$. If the regular grading is minimal and for every $h,h^{\prime}\in H$ and every $0\neq r_h\in R_h$ there exists $s_{h^{\prime}}\in R_{h^{\prime}}$ such that $r_{h}s_{h^{\prime}}\neq 0$ then $Z(R)=R_\epsilon$.
\end{proposition}
\textit{Proof.}
Since $\beta$ is a bicharacter it is clear that $\beta(\epsilon, h)=1$ for every $h\in H$. Condition (P2) of the above definition now implies that $R_{\epsilon}\subseteq Z(R)$. The center $Z(R)$ is a homogeneous subalgebra, hence if $R_{\epsilon} \subsetneqq Z(R)$ there exists $h\neq \epsilon$ and $0\neq z_h\in R_h$ such that $z_{h}\in Z(R)$. For $h^{\prime}\in H$ let $w_{h^{\prime}}$ be such that $z_hw_{h^{\prime}}\neq 0$. We have $w_{h^{\prime}}z_h=z_hw_{h^{\prime}}=\beta(h,h^{\prime})w_{h^{\prime}}z_h$. Hence $\beta(h,h^{\prime})=1$ for every $h^{\prime}\in H$. This is a contradiction because the regular grading is minimal.
\hfill $\Box$

\begin{remark}
None of the two hypotheses of the proposition previous can be removed. Indeed, if $R$ is an algebra with regular grading is not minimal, then there exist $\epsilon_H \neq h\in H$ such that $\beta(h,g)\neq 1$, for all $g\in H$ and so $Z(R)$ contains properly $R_\epsilon$.

We consider the algebra $R$ with presentation: $$\langle 1, f, e_1, e_2, \dots | fe_i=e_if=0, e_je_j=-e_je_i \rangle.$$ The subalgebra of $R$ generated by the $e_i$ is the Grassmann algebra $E$. The decomposition 
$$R_0=E_0\ \ \ \mbox{ and } \ \ \ R_1=E_1\oplus\langle f\rangle$$ is a minimal regular $\mathbb{Z}_2$-grading on $R$. It is clear that $f\in Z(A)$, hence $Z(R)$ contains properly $R_\epsilon$.
\end{remark}

\begin{lemma}\label{V1}
Let $R$ be an algebra with a regular grading by an abelian group $H$ such that $Z(R)=R_{\epsilon}$ and let $A$ be the algebra $M_n(R)$ with a $G$-grading. If $f=f(x_{1,g_1},\dots, x_{r,g_r})$ is a polynomial such that the result of every admissible substitution lies in $R_h$, for some $h\neq\epsilon_H$, then $f$ is a graded identity for $A$.
\end{lemma}
\textit{Proof.}
Let $(a_1,\dots, a_r)$ be an admissible substitution, Corollary \ref{CV} implies that for every $h^\prime\in H$ and every $m^{\prime}\in R_{h^\prime}$ the commutator $[f(a_1,\dots, a_r), m^{\prime}]$ lies in $R_{h}$. The group $H$ is abelian, hence it is clear that $f(a_1,\dots, a_r)\cdot m^{\prime}$ and $m^{\prime}\cdot f(a_1,\dots, a_r)$ lie in $R_{hh^{\prime}}$. Hence $[f(a_1,\dots, a_r), m^{\prime}]$ also lies in $V_{hh^\prime}$. If ${h^\prime}\neq \epsilon_H$ then $h\neq{hh^\prime}$ and this implies $V_{hh^\prime}\cap V_h=\{0\}$. Therefore we conclude that $[f(a_1,\dots, a_r), b]=0$. Hence $f(a_1,\dots, a_r)$ commutes with every $m^{\prime}$ in $R_{h^\prime}$ if $h^{\prime}\neq \epsilon$. Since $\beta(\epsilon, h)=1$ we conclude that $f(a_1,\dots, a_r)$ commutes with every element in $R_{\epsilon}$. This implies that $f(a_1,\dots, a_r)$ lies in $Z(R)=R_{\epsilon_H}$. In this case $f(a_1,\dots, a_r)$ lies in $R_h\cap R_{\epsilon_H}=\{0\}$, therefore $f(a_1,\dots, a_r)=0$. Hence $f$ is a graded identity for $A$.
\hfill $\Box$

Let $R$ be an algebra with a regular grading by the abelian group $H$ and let $\beta:H\times H \rightarrow F^{\times}$ denote its commutation function. If $m(x_{1},\dots, x_{s})$ is a monomial of degree $d_i$ in $x_{i}$ and $(r_1,\dots, r_s)$ is an $r$-tuple of homogeneous elements of $R$ then it follows from the condition (P2) of Definition \ref{def.rel} that $m(r_{1},\dots, r_{s})= \lambda r_1^{d_1}\cdots r_s^{d_s}$, for some non-zero scalar $\lambda$.

\begin{remark}\label{rr}
The scalar $\lambda$ in the previous equality depends only on the monomial $m$, the commutation function $\beta:H\times H \rightarrow F^{\times}$  and the $s$-tuple $\textbf{h}=(h_1,\dots, h_s)$ such that $a_i\in R_{h_i}$. This scalar will be denoted $\epsilon^R_{\textbf{h}, m}$ hence we have \[m(r_{1},\dots, r_{s})= \epsilon^R_{\textbf{h}, m} r_1^{d_1}\cdots r_s^{d_s}.\]
\end{remark}

\begin{definition} Let $f(x_{1,g_1},\dots,x_{n,g_n})$ be a multihomogeneous polynomial in the algebra $F\langle X_G\rangle$ and $\textbf{h}=(h_1,\dots, h_n)\in H^n$. We write $f=\sum_i \alpha_i m_i$ as a linear combination of monomials and define the polynomial $f_{\textbf{h}}$ as
$$f_{\textbf{h}}(x_{1,g_1},\dots,x_{n,g_n})=\sum_{i}\epsilon^R_{\textbf{h}, m_i}\alpha_i m_i.$$
\end{definition}

\begin{remark}\label{subs0}
If $f(x_{1,g_1},\dots,x_{r,g_r})$ is a multihomogeneous polynomial of degree $d_i$ in $x_{i,g_i}$, $a_i\in R_{h_i}$ and $B_i\in M_n(F)$ then $$f(a_1B_1,\dots,a_rB_r)=a_1^{d_1}\cdots a_r^{d_r}f_{\textbf{h}}(B_1,\dots, B_r),$$ where $\textbf{h}=(h_1,\dots,h_r)$.
\end{remark}

Next we consider an algebra $R$ an algebra with a regular $H$-grading that satisfies the following condition:

\begin{enumerate}
\item[(R1)] For any two monomials $m_1,m_2\in F\langle X_H \rangle$, in distinct variables, that are not $H$-graded identities for $R$ the product $m_1\cdot m_2$ is not a graded identity for $R$.
\end{enumerate}

The regular gradings in \cite{AO}, \cite{Regev&Bah} satisfy this property. An algebra $R$ satisfying (R1) satisfies the graded analogue of the definition of a verbally prime algebra. 

\begin{proposition}\label{grprime}
Let $R$ be an algebra with a regular $H$-grading that satisfies (R1). If $f,g \in F\langle X_H \rangle$ are non-identities in distinct variables for $R$ then $f\cdot g$ is a non-identity for $R$.
\end{proposition}
\textit{Proof.}
Since the field is infinite we may assume without loss of generality that $f,g$ are multihomogeneous. It follows from Remark \ref{rr} that $f,g$ are, modulo $Id_H(R)$, scalar multiples of monomials in distinct variables. Hence (R1) implies that $f\cdot g$ is a non-identity for $R$.
\hfill $\Box$

\begin{corollary}\label{ccenter}
Let $R$ be an algebra with a minimal regular $H$-grading that satisfies (R1). If $R^{\prime}=F\langle X_H \rangle / Id_H(R)$, with its canonical $H$-grading, then $Z(R^{\prime})=R^{\prime}_{\epsilon}$.
\end{corollary}
\textit{Proof.}
Let $h, h^{\prime} \in H$ and $r^{\prime}_h=f+Id_H(R)$ a non-zero element in $R_h^{\prime}$. Let $m$ be such that $x_{m,h^{\prime}}$ and $f$ are in distinct variables. Proposition \ref{grprime} implies that $f\cdot x_{m,h^{\prime}}$ does not lie in $Id_H$. Therefore $r^{\prime}_hs^{\prime}_{h^{\prime}}\neq 0$, where $s^{\prime}_{h^{\prime}}= x_{m,h^{\prime}}+Id_H(R)$. The result now follows from Proposition \ref{center}.
\hfill $\Box$

In the next propositon we remark that the algebra $R$ above is verbally prime. Over an algebraically closed field of characteristic zero there is a description of algebras that have a nondegenerate regular grading by a finite group in \cite[Corollary 42]{AO}.

\begin{proposition}
An algebra $R$ with a regular $H$-grading that satisfies (R1) must be verbally prime.
\end{proposition}
\textit{Proof.}
Let $f(x_1,\dots, x_t)$, $g(x_{t+1},\dots, x_s)$ be two non-identities, in distinct variables, for $R$. Proposition \ref{C} implies that there exists multihomogeneous polynomials $f^{\prime}(x_1,\dots, x_{t^{\prime}})\in \langle f \rangle_T$, $g^{\prime}(x_{t^{\prime}+1},\dots, x_{s^{\prime}})\in \langle f \rangle_T$ and homogeneous elements $r_1,\dots, r_{t^{\prime}},\dots, r_{s^{\prime}}$ such that $f^{\prime}(r_1,\dots, r_{t^{\prime}})\neq 0$ and $g(r_{t^{\prime}+1},\dots, r_{s^{\prime}})\neq 0$. Remark \ref{rr} and (R1) imply that $f^{\prime}(x_1,\dots, x_{t^{\prime}})\cdot g^{\prime}(x_{t^{\prime}+1},\dots, x_{s^{\prime}})$ is not an identity for $R$. Since $f^{\prime}\cdot g^{\prime}$ lies in $\langle f\cdot g \rangle_T$ we conclude that $f\cdot g$ is not an ordinary identity for $R$. 
\hfill $\Box$

In the next theorem we consider an algebra $A=M_n(R)$ with a $G$-grading such that:

\begin{enumerate}
\item[(Gr1)] $R$ is a homogeneous subalgebra and $R\subseteq A_{\epsilon}$;
\item[(Gr2)] $M_n(F)$ is a homogeneous subalgebra.
\end{enumerate}

Now we state our main result of the section.

\begin{theorem}\label{MnA}
Let $A$ be the algebra $M_n(R)$ with a $G$-grading that satisfies (Gr1), (Gr2). If $R$ has a regular grading, by an abelian group $H$, satisfying (R1) and the algebra $B=M_n(F)$ (with the inherited $G$-grading) satisfies the primeness property for graded central polynomials then $A$ also satisfies the primeness property.
\end{theorem}
\textit{Proof.}
We first prove that we may assume without loss of generality that the regular grading on $R$ satisfies (R1) and $Z(R)=R_{\epsilon}$. Let $R$ with a regular $H$-grading that satisfies (R1), then every coarsening of this grading satisfies (R1). Every regular grading on $R$ admits a coarsening by a homomorphic image of $H$ that is minimal, hence we may assume that $R$ has a minimal regular grading that satisfies (R1). For simplicity we still denote the grading group by $H$. The canonical $H$-grading on the algebra $R^{\prime}=F\langle X_H \rangle/Id_H(R)$ is a minimal regular grading that satisfies (R1). Corollary \ref{ccenter} implies that $Z(R^{\prime})=R^{\prime}_{\epsilon}$. The algebra $A^{\prime}=M_n(R^{\prime})$  with the $G$-grading such that $A^{\prime}_{g}=\{r^{\prime}M|r^{\prime}\in R^{\prime}, M\in B_g\}$ satisfies the same $G$-graded identities as $A$. Hence $A$ has the primeness property if and only if $A^{\prime}$ does. Henceforth we assume that $R$ has a regular grading by $H$ such that (R1) and $Z(R)=R_{\epsilon}$ hold.

Let $f(x_{1,g_1},\dots,x_{r,g_r})$ such that the product $f\cdot g$ is a central polynomial for $A$ for some $g(x_{{r+1},g_{r+1}},\dots,x_{s,g_s})$ in distinct variables than $f$. We will prove that $f$ is a graded central polynomial for $A$. We first claim that there exists $g^{\prime}(x_{1,g_1},\dots,x_{m,g_m})$ in $\langle g \rangle_{T_G}$, an $m$-tuple $\textbf{h}=(h_1,\dots, h_m)$ in $H^{m}$ and an $g_{\textbf{h}}^{\prime}$-admissible substitution $(A_1,\dots, A_m)$ such that $g_{\textbf{h}}^{\prime}(A_1,\dots, A_m)$ is an invertible matrix. Let $\mathcal{B}_R$ be a basis of $R$ of homogeneous elements in the regular $H$-grading and $\mathcal{B}_B$ a basis of $B$ of homogeneous matrices in the $G$-grading. It is clear that $\mathcal{B}_A=\{mM| m \in \mathcal{B}_R, M \in \mathcal{B}_B\}$ is a basis for $A$. The product $f\cdot g$ is not an identity, hence Proposition \ref{C} implies that there exists $f^{\prime}(x_{1,g_1},\dots, x_{m, g_m})$ in $\langle f \rangle_{T_G}$ and $g^{\prime}(x_{1,g_1},\dots, x_{m, g_m})$ in $\langle g \rangle_{T_G}$ and an admissible substitution $(a_1A_1,\dots, a_mA_m)$ by elements of $\mathcal{B}_A$ such that
\begin{equation}\label{notzero}
f^{\prime}(a_1A_1,\dots, a_mA_m)\cdot g^{\prime}(a_1A_1,\dots, a_mA_m)\neq 0.
\end{equation}
Proposition \ref{ll} implies that the element above lies in the centre of $A$. Let $h_i$ denote the degree of $a_i$ in the regular grading of $R$ and $\textbf{h}=(h_1,\dots, h_m)$. Using Remark \ref{subs0} we rewrite the above element as \[(a_1^{d_1}\cdots a_m^{d_m})\cdot (a_1^{e_1}\cdots a_m^{e_m}) f^{\prime}_{\textbf{h}}(A_1,\dots, A_m)\cdot g^{\prime}_{\textbf{h}}(A_1,\dots, A_m),\] where $d_i$ and $e_i$ denote the degree of $x_{i,g_i}$ in $f^{\prime}$ and $g^{\prime}$, respectively. Since this element is central the product $f^{\prime}_{\textbf{h}}(A_1,\dots, A_m)\cdot g^{\prime}_{\textbf{h}}(A_1,\dots, A_m)$ is a scalar matrix, different from $0$ because of (\ref{notzero}). Hence $ g^{\prime}_{\textbf{h}}(A_1,\dots, A_m)$ is an invertible matrix and the claim is proved.

Let $h^{-1}$ be the degree of $(a_1^{e_1}\cdots a_m^{e_m})$ in the regular grading of $R$. We claim that the result of every admissible substitution in the polynomial $f$ lies $R_{h}$. It is clear that $f$ is not an identity, hence this claim and Lemma \ref{V1} imply that $f$ is central. Let $f^{\prime \prime}(y_{1,g_1},\dots, y_{s, g_s})$ be a polynomial in $\langle f \rangle_{T_G}$ and $(b_1B_1,\dots, b_sB_s)$ an admissible substitution by elements of $\mathcal{B}_A$. Denote $\textbf{k}=(k_1,\dots, k_s)$ the $s$-tuple where $k_i$ is the homogeneous degree of $b_i$ in the regular grading of $R$. The product $f^{\prime \prime}(y_{1,g_1},\dots, y_{s, g_s})\cdot g^{\prime}(x_{1,g_1},\dots,x_{m,g_m})$ is a central polynomial. Hence the element
\begin{equation}\label{ec}
(b_1^{f_1}\cdots b_s^{f_s})\cdot (a_1^{e_1}\cdots a_m^{e_m}) f^{\prime\prime}_{\textbf{k}}(B_1,\dots, B_k)\cdot g^{\prime}_{\textbf{h}}(A_1,\dots, A_m),
\end{equation}
where $f_i$ is the degree of $y_{i,g_i}$ in $f^{\prime \prime}$, is central. We will prove that $f^{\prime\prime}(b_1B_1,\dots, b_sB_s)$ lies in $R_h$, this and Lemma \ref{C} imply the claim. We may assume that
$$ f^{\prime\prime}(b_1B_1,\dots, b_sB_s)=(b_1^{f_1}\cdots b_s^{f_s}) f^{\prime\prime}_{\textbf{k}}(B_1,\dots, B_k)\neq 0.$$ Since the product of two monomials $m_1,m_2 \in F\langle X_H \rangle$ that are not identities for $R$ is also not an identity, using (R1) we may assume without loss of generality that
$$(b_1^{f_1}\cdots b_s^{f_s})\cdot (a_1^{e_1}\cdots a_m^{e_m})\neq 0.$$
In this case $ f^{\prime\prime}_{\textbf{k}}(B_1,\dots, B_m)\cdot g^{\prime}_{\textbf{h}}(A_1,\dots, A_m)$ is a scalar matrix. The previous argument proves that the result of every admissible substitution in the polynomial $f^{\prime\prime}_{\textbf{k}}\cdot g^{\prime}_{\textbf{h}}$ is a scalar matrix. The matrix $ f^{\prime\prime}_{\textbf{k}}(B_1,\dots, B_m)$ is not $0$ and the matrix $ g^{\prime}_{\textbf{h}}(A_1,\dots, A_m)$ is invertible, hence  $f^{\prime\prime}_{\textbf{k}}\cdot g^{\prime}_{\textbf{h}}$ is not an identity.  Since $B$ satisfies the primeness property we conclude that $f^{\prime\prime}_{\textbf{k}}$ is a central polynomial for $B$. Therefore $f^{\prime\prime}_{\textbf{k}}(B_1,\dots, B_k)$ is a scalar matrix. Since $Z(R)=R_{\epsilon}$ it is clear that the centre of $A$ is $R_{\epsilon}$. The element in (\ref{ec}) is central, hence it follows that $(b_1^{f_1}\cdots b_s^{f_s})\cdot (a_1^{e_1}\cdots a_m^{e_m})$ lies in $R_{\epsilon}$. Therefore $(b_1^{f_1}\cdots b_s^{f_s})$ lies in $R_{h}$. Hence conclude that $f^{\prime\prime}(b_1B_1,\dots, b_sB_s)$ lies in $R_{h}$.
\hfill $\Box$

\begin{remark}
Over a field of characteristic zero one may reduce the above argument to multilinear polynomials instead of multihomogeneous ones. In this case the no extra hypothesis is needed on the regular grading on $R$ to prove the previous theorem. Indeed condition (R1) for multilinear monomials is a direct consequence of (P1) in Definition \ref{def.rel}.
\end{remark}

As an immediate consequence of the previous theorem we have the following corollaries.

\begin{corollary}\label{MnE}
The algebra $M_n(E)$ (with the trivial grading) satisfies the primeness property for central polynomials.
\end{corollary}
\textit{Proof.}
The canonical $\mathbb{Z}_2$-grading on $E$ is regular and satisfies the hypothesis of Theorem \ref{MnA}. Moreover, the algebra $M_n(F)$ with the trivial grading satisfies the primeness property (see \cite{Regev}).
\hfill $\Box$

The previous corollary generalizes to infinite fields the main result of \cite{D} which is valid for a field of characteristic zero.

\begin{corollary}\label{MnFcross}
Let $G$ be a group with no non-trivial homomorphism $G\rightarrow F^{\times}$. If the algebra $M_n(E)$ has a $G$-grading such that $E$ is a homogeneous subalgebra with the trivial grading and $M_n(F)$ a homogeneous subalgebra with the crossed product grading then $M_n(E)$ satisfies the primeness property for graded central polynomials.
\end{corollary}
\textit{Proof.}
The canonical $\mathbb{Z}_2$-grading on $E$ is regular and satisfies the hypothesis of Theorem \ref{MnA}. The result now follows from Theorems \ref{MnA} and \ref{maing}.
\hfill $\Box$

The algebra $R=M_m(F)$ has a regular grading that satisfies the hypothesis of Theorem \ref{MnA} (see \cite{Regev&Bah}). Since $M_{n}(R)\cong M_{mn}(F)$ one may produce other examples of gradings on matrix algebras with entries in the field $F$ that satisfy the primeness property for graded central polynomials using the main result of this section.

\section{Primeness property for the algebras $M_{a,b}(E)$}\label{smab}

The verbally prime algebras $M_n(F)$ and $M_n(E)$ have regular gradings satisfying (R1), (R2). Therefore the primeness property for these algebras may be obtained as a consequence of Theorem \ref{MnA} (with n=1). The remaining type of verbally prime algebra introduced by A. Kemer is $M_{a,b}(E)$. If $a\neq b$ no regular grading is known for this algebra. In this section we use the canonical $\mathbb{Z}_2$-grading and the results in Section \ref{s1} to prove that $M_{a,b}(E)$ satisfies the primeness property for ordinary central polynomials. 

\begin{definition}\label{envel}
Let $A$, $R$ be two $H$-graded algebras. We denote by $A\widehat{\otimes}R$ the $H$-graded
algebra such that $(A\widehat{\otimes}R)_{h}=A_{h}\otimes R_h$.
\end{definition}

The algebra $M_{a,b}(E)$ is obtained as $M_{a+b}(F)\widehat{\otimes}E$ where $M_{a+b}(F)$ has the elementary $\mathbb{Z}_2$-grading induced by the tuple $(0,\dots,0,1,\dots,1)$ with $a$ (resp. $b$) entries equal to $0$ (resp. $1$). If $R$ has a regular grading there is a natural relation between the graded identities of $A$ and $A\widehat{\otimes}R$ (see \cite{AO}, \cite{BD}).

Given a monomial $m(x_{1,h_1},\cdots, x_{k,h_k})$ we denote $\epsilon_m^{R}$ the scalar in Remark \ref{rr}.

\begin{definition}
Let $f(x_{1,h_1},\dots, x_{n,h_n})$ be a polynomial in $F\langle X_H \rangle$ and write $f=\sum \lambda_im_i$ as a linear combination of monomials $m_1,\dots, m_n$. We denote $f^{*}$ the polinomial $f^{*}=\sum \epsilon_{m_i}^{R}\lambda_im_i$.
\end{definition}

Over a field of characteristic zero we have the following result.

\begin{proposition}\cite[Lemma 27]{AO}
Let $F$ be a field of characteristic zero, $A$ be an $H$-graded algebra and $R$ an algebra with a regular $H$-grading. If $f$ is a multilinear polynomial then $f$ is a graded identity for $A$ if and only if $f^{*}$ is a graded identity for $A\widehat{\otimes}R$.
\end{proposition}

\begin{remark}\label{inv}
The "envelope operation" in Definition \ref{envel} is involutive. More precisely if $H$ is finite and $\widehat{R}=\widehat{\otimes}^{|H|-1}R$ then $\widehat{R}$ has a regular $H$-grading and $Id_H((A\widehat{\otimes}R)\widehat{\otimes}\widehat{R})=Id_{H}(A)$.This is proved in \cite[Theorem 6]{AO}.
\end{remark}

In this case we obtain the following theorem.

\begin{theorem}
Let $F$ be a field of characteristic zero, $A$ be a graded algebra graded by a finite abelian group $H$ and $R$ an algebra with a regular $H$-grading. The algebra $A$ has the primeness property for $H$-graded central polynomial if and only if $A\widehat{\otimes}R$ has the primeness property.
\end{theorem}
\textit{Proof.}
Let $f(x_{1,h_1},\dots, x_{r,h_r})$, $g(x_{r+1,h_{r+1}},\dots, x_{s,h_s})$ be multilinear polynomials in distinct variables. It is clear that $(f\cdot g)^{*}=f^{*}\cdot g^{*}$. Note also that $f\mapsto f^{*}$ is an invertible linear map in the subspace of multilinear polynomials in a fixed set of variables. Hence every multilinear polynomial is of the form $h^{*}$ for some multilinear $h$. We assume that $A$ has the primeness property. To prove that $A\widehat{\otimes}R$ we may to consider multilinear polynomials only, hence suffices to prove that if $f^{*}\cdot g^{*}$ is central then $f^{*}, g^{*}$ are central. The previous proposition implies that $f$ is central for $A$ if and only if $f^{*}$ is central for $A\widehat{\otimes}R$. Since $f^{*}\cdot g^{*}$ is central for $A\widehat{\otimes}R$ and $(f\cdot g)^{*}=f^{*}\cdot g^{*}$ we conclude that $f\cdot g$ is central for $A$. This implies that $f,g$ are central and therefore $f^{*},g^{*}$ are central. Hence $A\widehat{\otimes}R$ has the primeness property. 

Now assume that $A\widehat{\otimes}R$ has the primeness property, by Remark \ref{inv} the $H$-envelope operation is involutive. Hence we conclude by the above implication that $A$ has the primeness property.
\hfill $\Box$

\begin{remark}
It follows from Theorem \ref{maing} that $M_2(F)$ with its canonical elementary $\mathbb{Z}_2$-grading does not satisfy the primeness property, hence by the above result $M_{1,1}(E)$ with its canonical $\mathbb{Z}_2$-grading does not satisfy the property. We may verify this directly: $x_{1,\overline{1}}^2\cdot x_{2,\overline{1}}^2$ is a graded central polynomial, however $x_{1,\overline{1}}^2$, $x_{2,\overline{1}}^2$ are not graded central polynomials.
\end{remark}

We conclude this section by proving that $M_{a,b}(E)$ satisfies the primeness property for ordinary central polynomials.

\begin{lemma}\label{lemma}
Let $A=A_0\oplus A_1$ denote the algebra $M_{a,b}(E)$ with its canonical $\mathbb{Z}_2$-grading and let  $f(x_1,\dots, x_r)$ be a polynomial such that for some $g(x_{r+1},\dots, x_s)$, in distinct variables than $f$, the product  $f\cdot g$ is a central polynomial for $M_{a,b}(E)$. Then there exists an $\overline{i}$ in $\mathbb{Z}_2$ such that $f(a_1,\dots, a_r)$ lies in $A_{\overline{i}}$ for every $a_1,\dots, a_r$ in $M_{a,b}(E)$.
\end{lemma}

\textit{Proof}
Let $\textbf{\emph{B}}$ denote the canonical basis of $A$. We prove first that there exists multihomogeneous polynomials $f^{\prime}(x_1,\dots, x_m)$ in $\langle f \rangle_T$ and $g^{\prime}(x_{1}, \dots, x_m)$ in $\langle g \rangle_T$  and elements $b_{1},\dots, b_m$ in $\textbf{\emph{B}}$ such that $f^{\prime}(b_1,\dots, b_m)\cdot g^{\prime}(b_{1}, \dots, b_m)\neq 0$. Otherwise the polynomials of the form $f^{\prime}\cdot g^{\prime}$ satisfy assertion $(iii)$ of Proposition \ref{C} with $V=0$. Since these polynomials generate $\langle f\cdot g \rangle_T$ as a vector space Proposition \ref{C} implies that $f\cdot g$ is an identity for $A$, which is a contradiction. The elements $f^{\prime}(b_1,\dots, b_m)$ and $g^{\prime}(b_{1}, \dots, b_m)$
have the form $mM$ and $rR$, where $M$ and $R$ are matrices in $M_n(F)$ and $m$ and $r$ are monomials in $E$. Proposition \ref{ll} implies that the non-zero product $(mM)(rR)$ lies in $Z(A)$. Hence $MR$ is a non-zero scalar matrix and therefore $R$ is an invertible matrix. Let $\overline{i}$ be the element of $\mathbb{Z}_2$ such that $g^{\prime}(b_{1}, \dots, b_m)=rR$ lies in $A_{\overline{i}}$. For a  multihomogeneous polynomial $f^{\prime \prime}(x_{i_1},\dots, x_{i_s})$ in $\langle f \rangle_T$ and $b_{i_1},\dots, b_{i_s}$ in $\textbf{\emph{B}}$ we write $f^{\prime \prime}(b_{i_1},\dots, b_{i_s})=qQ$, where $q$ is a monomial in $E$ and $Q$ a matrix in $M_n(F)$. If $qQ=0$ then it lies in $A_{\overline{i}}$. Now we suppose that $qQ\neq 0$, we may assume without loss of generality that $qr\neq 0$. Since $R$ is invertible and $Q\neq 0$ the product $QR$ is not $0$ and therefore $(qQ)(rR)\neq 0$. Proposition \ref{ll} implies that $f^{\prime \prime}(b_{i_1},\dots, b_{i_s})\cdot g^{\prime}(b_{1}, \dots, b_m)$ is a central polynomial for $A$. Therefore $(qQ)(rR)$ lies in $Z(A)\subseteq A_0$. Note that $qQ$ lies in $A_{\overline{j}}$ for some $\overline{j}$ in $\mathbb{Z}_2$. Hence $(qQ)(rR)$ also lies in $A_{\overline{j} + \overline{i}}$. Since $(qQ)(rR)\neq 0$ this implies that $\overline{j} + \overline{i}=0$. Therefore $qQ$ lies in $A_{\overline{i}}$. Proposition \ref{C} with $V=A_{\overline{i}}$ yields the result.
\hfill $\Box$

\begin{theorem}\label{main}
Let $F$ be an infinite field of characteristic different from $2$. The $F$-algebra $M_{a,b}(E)$ has the primeness property for central polynomials.
\end{theorem}

\textit{Proof}
Let $f$ and $g$ be polynomials in distinct variables such that  $f\cdot g$ is a central polynomial for $A=M_{a,b}(E)$. It follows from the previous proposition that there exists an $\overline{i}$ in $\mathbb{Z}_2$ such that $f(a_1,\dots, a_r)$ lies in $A_{\overline{i}}$ for every $a_1,\dots, a_r$ in $M_{a,b}(E)$. The previous proposition and Corollary \ref{CV} imply that for every $b$ in $A$ the commutator $[f(a_1,\dots, a_r), b]$ also lies in $A_{\overline{i}}$. If $b$ lies in $A_1$ then the products $f(a_1,\dots, a_r)\cdot b$ and $b\cdot f(a_1,\dots, a_r)$ lie on $A_{\overline{i}+1}$. Hence the commutator $[f(a_1,\dots, a_r), b]$ lies on $A_{\overline{i}+1}$. We conclude that $[f(a_1,\dots, a_r), b]$ lies on $A_{\overline{i}}\cap A_{\overline{i}+1}$, hence $[f(a_1,\dots, a_r), b]=0$. To prove that $f(a_1,\dots, a_r)$ commutes with every element in $A_0$ it suffices to prove that it commutes with elements $b=e_{i_1}\cdots e_{i_{2k}}E_{ii}$, for $k>0$. We choose $j$ such that $b_1=e_{i_1}E_{ij}$ and $b_2=e_2\cdots e_{i_{2k}}E_{ji}$ lie in $A_1$. Since $b=b_1\cdot b_2$ we conclude that $f(a_1,\dots, a_r)$ commutes with $b=e_{i_1}\cdots e_{i_{2k}}E_{ii}$. Therefore the element $f(a_1,\dots, a_r)$ is central.  The proof that $g$ is a central polynomial is analogous.
\hfill $\Box$

\begin{remark}
Over a field of characteristic zero  the only verbally prime algebras are $M_n(F)$, $M_n(E)$ and $M_{a,b}(E)$ (see \cite{Kemer2}) and we conclude that every verbally prime algebra satisfies the primeness property for central polynomials.
\end{remark}

%\begin{flushleft}
%\textbf{Acknowledgements}
%\end{flushleft}

\end{document}